\documentclass[12pt,leqno]{article}
\usepackage[pdftex]{graphicx}
\usepackage{amsfonts,amssymb,amsmath}
\usepackage{ccaption,mathrsfs,amsmath,amssymb,amsthm}
\usepackage{amsmath, amscd, amsfonts,  tikz, times, wrapfig}
\usepackage{caption}

\newcounter{conjecture}
\newcounter{remark}

\newcommand{\eqnsection}{
    \renewcommand{\theequation}{\thesection.\arabic{equation}}
    \makeatletter
    \csname @addtoreset\endcsname{equation}{section}
    \makeatother}
\newtheorem{theorem}{Theorem}
\newtheorem{defn}{Definition}

\newtheorem{prop}{Proposition}

\newcommand{\dd}{\delta}

\newcommand{\lar}{\longrightarrow}
\newcommand{\eps}{\varepsilon}

\newcommand{\CC}{\mathbb{C}}

\def \be{\begin{equation}}
\def \ee{\end{equation}}
\def \bt{\begin{theorem}}
\def \et{\end{theorem}}
\def \bea{\begin{eqnarray}}
\def \eea{\end{eqnarray}}
\def \bas{\begin{eqnarray*}}
\def \eas{\end{eqnarray*}}

\newcommand {\rrr}[1]{(\ref{#1})}

\def \la{\lambda}

\def \Om{\Omega}

\def \si{\sigma}

\def \th{\theta}

\def \ff{\infty}

\def \DD{{\mathbb D}}

\def \HH{{\mathbb H}}

\def \RR{{\mathbb R}}

\def \ZZ{{\mathbb Z}}

\def \({\left(}
\def \){\right)}

\def \vski{\vspace{12pt}}

\newcommand{\bsh}{\backslash}

\def \bc{\begin{center} }
\def \ec{\end{center} }
\def \bs{\begin{slide} }
\def \es{\end{slide} }

\def\square{{\vcenter{\vbox{\hrule height.3pt
         \hbox{\vrule width.3pt height5pt \kern5pt
            \vrule width.3pt}
         \hrule height.3pt}}}}
\def\qed{{\hfill $\Box$ \bigskip}}

\newcounter{ccases}
\eqnsection
\begin{document}

\title{On the planar Brownian Green's function for stopping times.}

\author{
\begin{tabular}{c}
\textit{Greg Markowsky} \\
Monash University \\
Victoria 3800, Australia \\
gmarkowsky@gmail.com
\end{tabular}}

\bibliographystyle{amsplain}

\maketitle \eqnsection \setlength{\unitlength}{2mm}

\begin{abstract}
\noindent It has been known for some time that the Green's function of a planar domain can be defined in terms of the exit time of Brownian motion, and this definition has been extended to stopping times more general than exit times. In this paper, we extend the notion of conformal invariance of Green's function to analytic functions which are not injective, and use this extension to calculate the Green's function for a stopping time defined by the winding of Brownian motion. These considerations lead to a new proof of the Riemann mapping theorem. We also show how this invariance can be used to deduce several identities, including the standard infinite product representations of several trigonometric functions.
\vski


\vski

\end{abstract}

\section{Introduction}

In the field of analysis, Green's function $G(x,y)$ on regions of $\RR^n$ is formally defined to be the solution of $LG(x,y) = \dd(y-x)$, where $L$ is a linear differential operator. In complex analysis, where the Laplacian is the differential operator of most importance, for a given domain $\Om \subseteq \CC$ and $z \in \Om$ the Green's function of the Laplacian is generally defined by the following.

\begin{defn} \label{anal}

The Green's function $G_\Om(z,w)$ on a domain $\Om$ is a function in $w$ on $\Om \bsh \{z\}$ satisfying the following properties.

\begin{itemize} \label{analdef}

\item[(i)] $G_\Om(z,w)$ is harmonic and positive on $\Om \bsh \{z\}$.

\item[(ii)] $G_\Om(z,w) \lar 0$ as $w \lar \dd \Om$ (Note that the boundary is to be taken in the Riemann sphere, so that if $\Om$ is unbounded then $\ff \in \Om$).

\item[(iii)] $G_\Om(z,w) + \frac{1}{\pi} \ln |w-z|$ extends to be continuous (and therefore harmonic) at $w=z$.

\end{itemize}

\end{defn}

In fact, as defined here, it can be shown by standard analytic techniques that $G_\Om$ satisfies $LG_\Om(z,w) = \frac{1}{\pi}\dd(w-z)$. Note that the normalization here has been chosen to align with the probabilistic considerations to follow; removing the multiplicative constant $\pi^{-1}$ from $(iii)$ would result in a solution to $LG_\Om(z,w) = \dd(w-z)$. Not every domain has a Green's function as defined above, as for instance it can be shown that no such function can exist on the punctured disk $\DD^\times = \{0 < |z| < 1\}$, or more generally on a domain with isolated singularities. The Green's function is of tremendous importance in analysis on $\RR^n$, including complex analysis, and the question of what domains possess a Green's function has been keenly studied by analysts over the years.

\vski

On the other hand, the term "Green's function" has entered the vocabulary of probabilists in a way that may seem initially unrelated, namely as a measure the expected number of times that a discrete process visits a point, or the expected amount of time that a continuous process spends at a point. The reason that these two different notions have garnered the same name was discovered by Hunt in 1956 (\cite{hunt}), who showed that in many cases in $\RR^n$ these two notions coincide, with $L$ the Laplacian and the process in question Brownian motion. Before discussing this fact further, let us examine the probabilistic notion of a Green's function as pertains to Brownian motion in more detail. 

\vski

We will let a Brownian motion $B_t$ run until a stopping time $\tau$, and will let $B_t = \Delta$ for $t \geq \tau$, where $\Delta$ is a so-called "cemetery point" outside of $\hat \CC$. The point is that $B_t$ should no longer be in the plane for $t \geq \tau$. Let $\rho_t^\tau(z,w)$ be the probability density function at point $w$ and time $t$ of this killed Brownian motion. We then can calculate formally, for any measurable function $f$ on $\Om$,

\begin{equation} \label{exploctime}
\begin{split}
E_z \int_{0}^{\tau} f(B_s) ds & =  \int_{0}^{\tau} E_z[f(B_s)] ds \\
& = \int_{0}^{\ff} \int_{\CC} f(w) \rho_s^\tau(z,w) dA(w) ds \\
& = \int_{\CC} \Big( \int_{0}^{\ff} \rho_s^\tau(z,w) ds \Big) f(w) dA(w).
\end{split}
\end{equation}

This leads one to the consideration of the function

\begin{equation} \label{probdef}
G_\tau(z,w) := \int_{0}^{\ff} \rho_s^\tau(z,w) ds.
\end{equation}

This definition appears in a number of places, including \cite{hunt}, \cite{durBM}, \cite{bass}, and \cite{mortper}. Most often the definition has been examined with $\tau$ being the exit time from a domain, but as is noted in \cite{mortper} and \cite{bass} more general stopping times are allowable. Going forward, we will always use the notation $T_\Om$ to denote the exit time of a domain $\Om$; that is, $T_\Om = \inf\{t \geq 0 : B_t \in \Om^c\}$. In the case that $\tau = T_\Om$ for some domain $\Om$, we will simplify notation by writing $G_{\Om}(z,w) := G_{T_\Om}(z,w)$. Note that we are using the same notation here as in Definition \ref{analdef}, and similarly we will refer to both concepts as "Green's functions". In order to distinguish between them where necessary, we will call a function satisfying the conditions in Definition \ref{analdef} the {\it analyst's Green's function}, and we will refer to \rrr{probdef} as the {\it probabilist's Green's function}. The analyst's Green's function is well known to be conformally invariant: if $f$ is a conformal map on $\Om$, then $G_\Om(z,w) = G_{f(\Om)}(f(z),f(w))$; this is virtually immediate from the definition, as a harmonic function composed with a conformal one is again analytic, and conditions $(i)-(iii)$ follow easily. The probabilist's Green's function is conformally invariant as well, which is a simple consequence of the conformal invariance of Brownian motion, as we will show later in this section. That the two functions coincide on any simply connected domain can therefore be proved by (a) checking that they coincide on some domain, such as the unit disk or upper half-plane and (b) invoking the Riemann mapping theorem, which states that any two simply domains (excluding $\CC$ itself) are conformally equivalent; this argument would be certainly known to many researchers familiar with these topics. It would seem that this leaves little more to say about the probabilist's Green's function, since the analyst's Green's function is so well understood. However, the primary goal of this paper is to ignore this argument, and show the equivalence of the two notions of Green's functions on simply connected domains using only properties of Brownian motion and a few simple analytic functions. A consequence of this is an intuitive proof of the Riemann mapping theorem.

\vski

Our first order of business is showing that the probabilist's Green's function is conformally invariant, and subsequently providing an extension of this statement to analytic functions which are not necessarily injective. The main tool necessary is the following theorem due to L\'evy (see \cite{durBM} or \cite{mortper} for a proof).

\begin{theorem} \label{holinv}
Let $f$ be analytic and nonconstant on a domain $\Omega$, and let $a \in \Omega$.  Let $B_t$ be a Brownian motion started at $a$ and stopped at a stopping time $\tau$ such that the set of Brownian paths $\{B_t: 0 \leq t \leq \tau\}$ lie within the closure of $\Om$ a.s. Set

\begin{equation} \label{tokyo}
\sigma_t = \int_{0}^{t \wedge \tau} |f'(B_s)|^2 ds.
\end{equation}

Then there is a Brownian motion $\hat B$ starting at $f(a)$ and stopped at the stopping time $\sigma_{\tau}$ such that $f(B_t) = \hat B_{\si_t}$.
\end{theorem}

As was alluded to earlier, the following proposition and its proof must certainly be known, but it seems difficult to locate a simple statement in the literature (perhaps for wont of applications). For the benefit of the reader, we therefore provide a quick proof.

\begin{prop} \label{}
The probabilist's Green's function is conformally invariant. That is, if $f$ is a conformal map on $\Om$, then $G_\Om(z,w) = G_{f(\Om)}(f(z),f(w))$.
\end{prop}

{\bf Proof:} Define a measure $\mu_\Om$ on Borel subsets of $\Om$ by

\begin{equation} \label{}
\mu_\Om^z(A) = E_z\int_{0}^{T_\Om} 1_A(B_t)dt.
\end{equation}

Standard arguments now show that $\mu_\Om^z$ has a density equal to $G_\Om(z,w)$. If $D(w,\dd)$ denotes the disk of radius $\dd$ centered at $w$ and $\la$ denotes Lebesgue measure in the plane, then we have

\begin{equation} \label{}
G_\Om(z,w) = \lim_{\dd \searrow 0} \frac{\mu^z_\Om(D(w,\dd))}{\la(D(w,\dd))}; \qquad G_{f(\Om)}(f(z),f(w)) = \lim_{\dd \searrow 0} \frac{\mu^{f(z)}_{f(\Om)}(f(D(w,\dd)))}{\la(f(D(w,\dd)))}.
\end{equation}

Let $\eps > 0$ be given. Conformality implies that $f'(w) \neq 0$, and that for sufficiently small $\dd$ we have $D(f(w),|f'(w)|\dd(1-\eps)) \subseteq f(D(w,\dd))) \subseteq D(f(w),|f'(w)|\dd(1+\eps))$, so that $\frac{\la(f(D(w,\dd)))}{\la(D(w,\dd))} \in (|f'(w)|^2(1-\eps)^2,|f'(w)|^2(1+\eps)^2)$. On the other hand, \rrr{tokyo} shows that the scaling factor for time is $|f'|^2$ as well; that is, we can choose $\dd$ sufficiently small so that $\frac{|f'(w')|}{|f'(w)|} \in (1-\eps,1+\eps)$ for all $w' \in D(w,\dd)$, and then we will have $\frac{\mu_\Om(f(D(w,\dd)))}{\mu_\Om(D(w,\dd))} \in (|f'(w)|^2(1-\eps),|f'(w)|^2(1+\eps))$. Combining these two estimates and letting $\eps \lar 0$ shows that $\lim_{\dd \searrow 0} \frac{\mu_\Om(D(w,\dd))}{\la(D(w,\dd))} = \lim_{\dd \searrow 0} \frac{\mu_{f(\Om)}(f(D(w,\dd)))}{\la(f(D(w,\dd)))}$, and the result follows. \qed

It is important to note, however, that L\'evy's theorem in fact does not require maps to be injective, permitting general nonconstant analytic functions as well. This allows us to extend the previous proposition (with essentially the same proof) as follows.

\begin{prop} \label{chimass}
Let $\Om$ be a domain, and suppose $f$ is a function analytic on $\Om$. Let $B_t$ be a Brownian motion starting at $a$, and $\tau$ a stopping time such that the set of Brownian paths $\{B_t: 0 \leq t \leq \tau\}$ lie within $U$ a.s. Let $\hat \tau = \si_\tau$, where $\si$ is defined by \rrr{tokyo}. Then

\begin{equation} \label{cherry2}
G_{\hat \tau} (f(z),w) = \sum_{w' \in f^{-1}(\{w\})} n(f,w') G_\tau (z,w'),
\end{equation}

where $n(f,w')$ is the order of the zero of $f-w$ at $w'$.
\end{prop}

The proof of the previous proposition applies at all points which are not images under $f$ of a point $w'$ at which $f'(w')=0$ (the {\it critical values} of $f$), with the only difference being that a point $w$ with multiple preimages will accumulate mass in the projected Green's function corresponding to mass accumulated at each of the preimages by the initial Brownian motion. The $n(f,w')$ term is not necessary in order to calculate a density of $\mu_{\hat \tau}$, as in fact $n(f,w') = 1$ except on the zero set of $f'$, which is a discrete set and therefore of Lebesque measure 0; the term is required merely to make the right side of \rrr{cherry2} continuous in $w$ when $G_\tau$ is as well (such as when $\tau$ is an exit time of a domain), since if $w$ is a critical value of $f$ then points near $w$ will have as many preimages as does $w$ only when multiplicities of preimages are counted. An illustrative example of this is given in the next section.

\vski


%


In the next two sections, we will show how the Riemann mapping theorem can be proved using the probabilist's Green's function. The role of the next section is simply to calculate the function for a few simple stopping times that we will need for the proof, which will be given in the subsequent section. In the section following that, we show how basic properties of Brownian motion and judicious choices of domains can be used to prove some nontrivial identities, specifically several infinite product representations for trigonometric functions, some standard and some less so. A final section contains a few concluding remarks.

\section{A few simple examples}

In this section we calculate a few examples that are required for the ensuing sections. It should be emphasized here that we make use only of elementary properties of Brownian motion and a few standard analytic functions. We note first the two-dimensional Gaussian density $\rho_t^\CC (z,w) = \frac{1}{2\pi t} e^{-|z-w|^2/(2t)}$; it is evident that $G_{\CC}(z,w) = \ff$ for all $z,w$.

\vski

Let $\HH = \{Im(z)>0\}$. The density $\rho_t^\HH (z,w)$ must capture the probability of Brownian paths near $w$ at time $t$, but only those which have not previously intersected $\RR$. By the reflection principle for Brownian motion (see \cite{fima} or \cite{mortper}), the processes $B_t$ and

\be \hat B_t = \left \{ \begin{array}{ll}
B_t & \qquad  \mbox{if } t \leq T_\HH  \\
\overline B_t & \qquad \mbox{if } t > T_\HH \;,
\end{array} \right. \ee

have the same law. But $B_t$ is near $w$ and $t > T_\HH$ precisely when $\hat B_t$ is near $\bar w$ with $t > T_\HH$, and this occurs precisely when $\hat B_t$ is near $\bar w$, since the Brownian motion cannot travel from $z$ to $\bar w$ without first crossing $\RR$. We conclude that $\rho_t^\HH (z,w) = \rho_t^\CC (z,w) - \rho_t^\CC (z,\bar w) = \frac{1}{2\pi t} (e^{-|z-w|^2/(2t)} - e^{-|z-\bar w|^2/(2t)})$. We therefore have

\begin{equation} \label{hp}
\begin{split}
G_{\HH}(z,w) & = \int_{0}^{\ff} \frac{1}{2\pi t} (e^{\frac{-|z-w|^2}{2t}} - e^{\frac{-|z-\bar w|^2}{2t}}) dt \\
& = \frac{1}{2\pi} \int_{0}^{\ff} \frac{1}{t} \int_{|z-w|}^{|z-\bar w|} \frac{a}{t}e^{\frac{-a^2}{2t}} da dt \\
& = \frac{1}{2\pi} \int_{|z-w|}^{|z-\bar w|} \int_{0}^{\ff} \frac{a}{t^2} e^{\frac{-a^2}{2t}} dt da \\
& = \frac{1}{\pi} \int_{|z-w|}^{|z-\bar w|} \frac{1}{a} \int_{0}^{\ff} e^{-u} du da \\
& = \frac{1}{\pi} \ln \Big( \frac{|z-\bar w|}{|z-w|}\Big).
\end{split}
\end{equation}

$\phi(z) = -i \Big(\frac{z-1}{z+1}\Big)$ maps $\DD$ conformally onto $\HH$, sending $0$ to $i$. Thus,

\begin{equation} \label{}
\begin{split}
G_{\DD}(0,z) & =G_{\HH}(i,\phi(z)) = \frac{1}{\pi} \ln \Big( \frac{|i-i(\frac{\bar z -1}{\bar z + 1})|}{|i+i(\frac{z -1}{z + 1})|}\Big) \\
& = \frac{1}{\pi} \ln \Big| \frac{2z+2}{2|z|^2 + 2z}\Big| = \frac{1}{\pi} \ln \frac{1}{|z|}.
\end{split}
\end{equation}

To calculate $G_{\DD}(a,w)$ for arbitrary $a \in \DD$, let $\phi(z) = \frac{z-a}{1-\bar a z}$. It is a standard fact that $\phi$ is conformal map from $\DD$ onto itself. Thus,

\begin{equation} \label{lamb}
G_{\DD}(a,w) = G_{\DD}(\phi(a),\phi(w)) = \frac{1}{\pi} \ln \frac{|1-\bar a w|}{|w-a|}.
\end{equation}

As an aside, let us at this point return to \rrr{cherry2} and consider the need for including the multiplicity of preimages in the formula. For any nonzero $a, w \in \DD$, \rrr{lamb} gives $G_{\DD}(a^2,w^2) = \frac{1}{\pi} \ln \frac{|1-\bar a^2 w^2|}{|w^2-a^2|}$; on the other hand, the map $z \lar z^2$ takes $\DD$ to itself, and \rrr{cherry2} therefore implies (for $w \neq 0$)

\begin{equation} \label{}
G_{\DD}(a^2,w^2) = \frac{1}{\pi} \ln \frac{|1-\bar a w|}{|w-a|} + \frac{1}{\pi} \ln \frac{|1+\bar a w|}{|w+a|},
\end{equation}

since the preimages of $w^2$ under this map are $\pm w$. Clearly these values agree, however the point $w=0$ has only one primage under the map, and without taking into account the multiplicity of this preimage an application of \rrr{cherry2} in the same manner would lead to the contradiction $\ln \frac{1}{|a^2|} = \ln \frac{1}{|a|}$. Since the multiplicity of the preimage at $0$ is 2, however, the correct value is returned by \rrr{cherry2}.

\vski

Let us now calculate the Green's function of a stopping time defined by the winding of Brownian motion. If our Brownian motion $B_t$ starts at 1 then it will never hit 0 a.s., and we can then define $arg(B_t)$ as a continuous process (with $arg(B_0) = 0$). For any posi-

\begin{wrapfigure}{r}{8cm}
\vspace{-.7in}
\includegraphics[width=8cm]{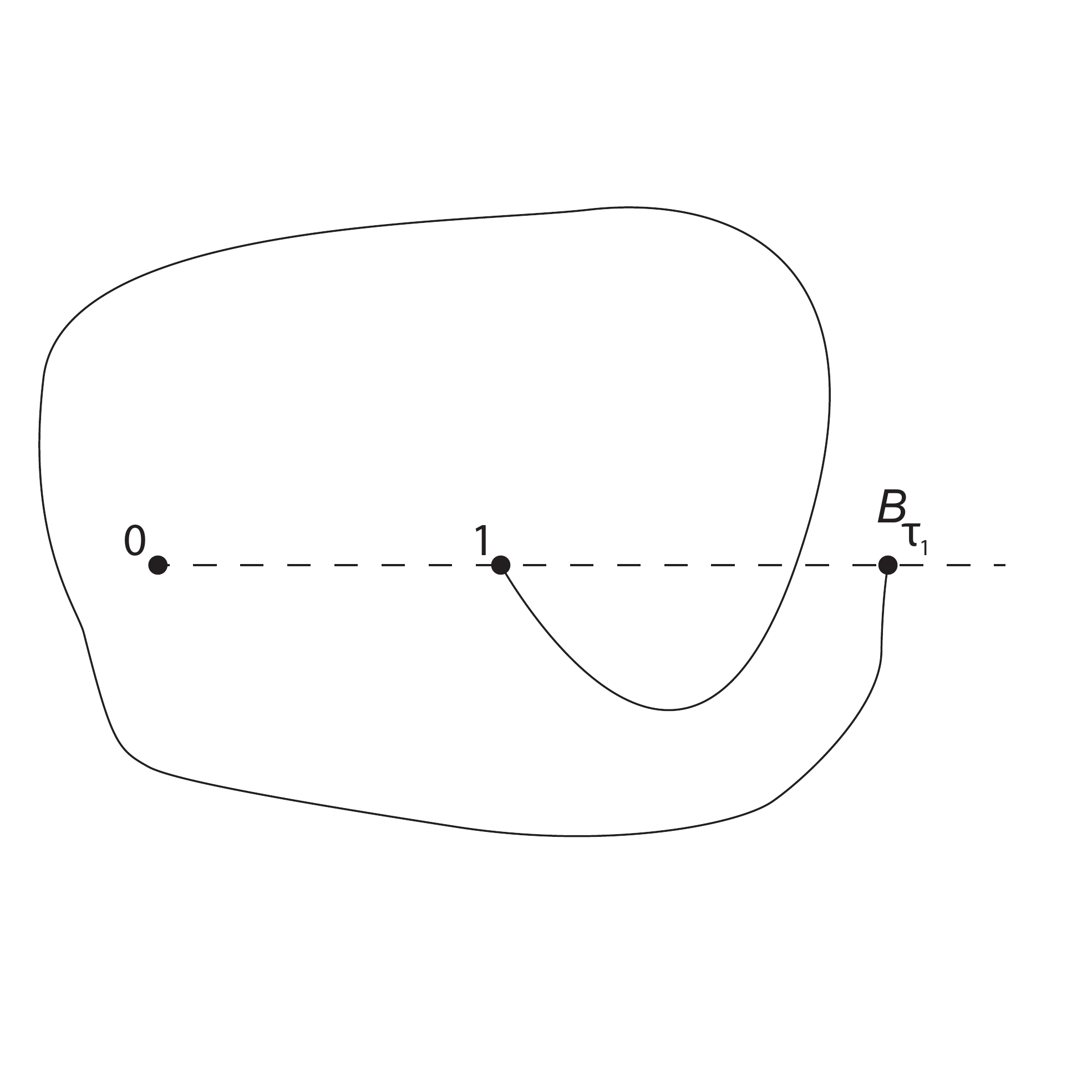}
\vspace{-.7in}
\caption{}
\end{wrapfigure}

tive integer $n$ we then define $\tau_n = \inf\{t: arg(B_t) = \pm 2 \pi n\}$. In other words, $\tau_n$ is the first time at which the Brownian motion has wound around the origin in either direction $n$ times. The image to the right gives an intuitive illustration of a Brownian path up until time $\tau_1$. If we let $U = \{Re(z)>0\}$ denote the right half-plane, we see that the function $z \lar z^{4n}$ transforms a Brownian motion starting at $1$ and stopped at $T_U$ into one starting at 1 and stopped at $\tau_n$. We may therefore calculate $G_{\tau_n}$ by using Proposition \ref{chimass}. For any point $re^{i\th}$, with $r >0$ and $\th \in (-\pi, \pi]$, we have

\begin{equation} \label{}
\nonumber G_{\tau_n}(1,re^{i\th}) = \sum_{\substack{k \in \ZZ \\ (\frac{\th}{4n}+ \frac{2\pi k}{4n}) \in [-\frac{\pi}{2},\frac{\pi}{2}]}} G_{U}(1,r^{\frac{1}{4n}} e^{i(\frac{\th}{4n}+ \frac{2\pi k}{4n})}).
\end{equation}

The explicit formula is now easy to produce if desired, using the fact that reflection across the $y$-axis is given by $z \lar -\bar z$, so that \rrr{hp} gives $G_U(z,w) = \frac{1}{\pi} \ln \frac{|z+\bar w|}{|z-w|}$; however, our purpose in examining this stopping time is the proof of the Riemann mapping theorem in the next section, and for that only a few basic properties of $G_{\tau_n}$ are required. To begin with, since $G_{\tau_n}$ is expressed as the sum of a finite number of finite terms, it is finite for all $re^{i\th} \neq 1$. Furthermore, $G_{\tau_n}(1,re^{i\th}) + \frac{1}{\pi} \ln |1-re^{i\th}|$ can be extended continuously at $re^{i\th} = 1$, since the only term in the sum which has a singularity corresponds to $k=0$ and has singular part $-\frac{1}{\pi} \ln |1-r^{\frac{1}{4n}} e^{i\frac{\th}{4n}}|$, and we have

\begin{equation} \label{}
\lim_{re^{i\th} \lar 1} \frac{1}{\pi} \ln |1-re^{i\th}|-\frac{1}{\pi} \ln |1-r^{\frac{1}{4n}} e^{i\frac{\th}{4n}}| = \lim_{\zeta \lar 1} \ln \frac{1-\zeta^{4n}}{1-\zeta} = \ln 4n.
\end{equation}

Finally, if $r \lar 0$ then all of the preimages of $re^{i\th}$ approach 0 as well, and we conclude that $G_{\tau_n}(1,re^{i\th}) \lar 0$; an analogous argument shows that $G_{\tau_n}(1,re^{i\th}) \lar 0$ whenever $r \lar \ff$.

\section{Proof of the Riemann mapping theorem}

A simply connected domain is a domain in which any closed curve is homotopic to a point. The Riemann mapping theorem is as follows.

\vski

{\bf Riemann mapping theorem} If $\Om \subsetneq \CC$ is a simply connected domain, and $a \in \Om$, then there is a conformal map $f$ from $\Om$ onto $\DD$, with $f(a)=0$.

\vski

The existence of such a map is equivalent to the existence of the analyst's Green's function on $\Om$: under the assumption of the existence of $f$, it can be shown easily that $\ln|f(w)|$ satisfies $(i) - (iii)$ in Definition \ref{anal}, while a bit more work shows that if we have $G_\Om$ then a harmonic conjugate $H(w)$ can be found on $\Om$ so that $f(w) = e^{-\pi(G_\Om(a,w)+ i H(w))}$ maps $\Om$ conformally onto $\DD$. This argument is spelled out in detail in a number of places, including the reference \cite{walshhistory}, which contains an interesting historical account as well. As may be seen by consulting that reference, the classical proofs of the Riemann mapping theorem proceeded by proving the existence of the analyst's Green's function for any simply connected domain, and this is how we will proceed as well, by showing that the probabilist's Green's function satisfies the conditions of Definition \ref{anal}.

\vski

To be precise, we take the probabilist's Green's function $G_\Om$, and must show that $G_\Om$ exists (i.e. is finite), is harmonic, has a logarithmic singularity at $a$, and is 0 on the boundary. The intuitive idea is very simple: if a Brownian motion winds enough times around a boundary point, then simple connectivity will imply that the Brownian motion has exited the domain, and we can bound $G_\Om$ by $G_{\tau_n}$ as defined in the previous section. We will make free use of the conformal invariance of $G_\Om$, and will where convenient change coordinates by a M\"obius transformation (also known as a linear fractional transformation, see \cite{ahl}), which allows us to take any three points in the sphere to any other three points. We will also make use of a monotonicity property of the probabilist's Green's function, namely that if $\tau_1, \tau_2$ are two stopping times with $\tau_1 \leq \tau_2$ a.s., then we have $G_{\tau_1}(a,w) \leq G_{\tau_2}(a,w)$ for all $w$; this is evident from the fact that $\rho_t^{\tau_1}(a,w) \leq \rho_t^{\tau_2}(a,w)$.

\vski

Given $\Om \subsetneq \CC$ simply connected and $a \in \Om$, choose a line passing through $a$ which intersects at least one point on $\Om^c$ other than $\ff$; this is always possible since $\Om$ is not all of $\CC$. Take $p_1,p_2$ to be the points on the line closest to $a$ on either side; note that one of $p_1, p_2$ may be $\ff$, but not both. After a M\"obius change of coordinates, we may assume $a=1, p_1=0, p_2=\ff$, and this implies that the entire ray $\RR^+ = (0,+\ff)$ lies within $\Om$. We now claim that $T_\Om < \tau_1$, where $\tau_1$ was defined in the previous section. To see why this must be so, we note that if we connect the Brownian path $\{B_t:0 \leq t \leq \tau_1\}$ to the point $1$ by adding the line segment $[1,B_{\tau_1}]$, we will have created a closed curve which is not homotopic to a point in $\CC \bsh \{0\}$. If the curve remains entirely within $\Om$, then we have contradicted simple connectivity.

\begin{wrapfigure}{r}{8cm}
\vspace{-.7in}
\includegraphics[width=8cm]{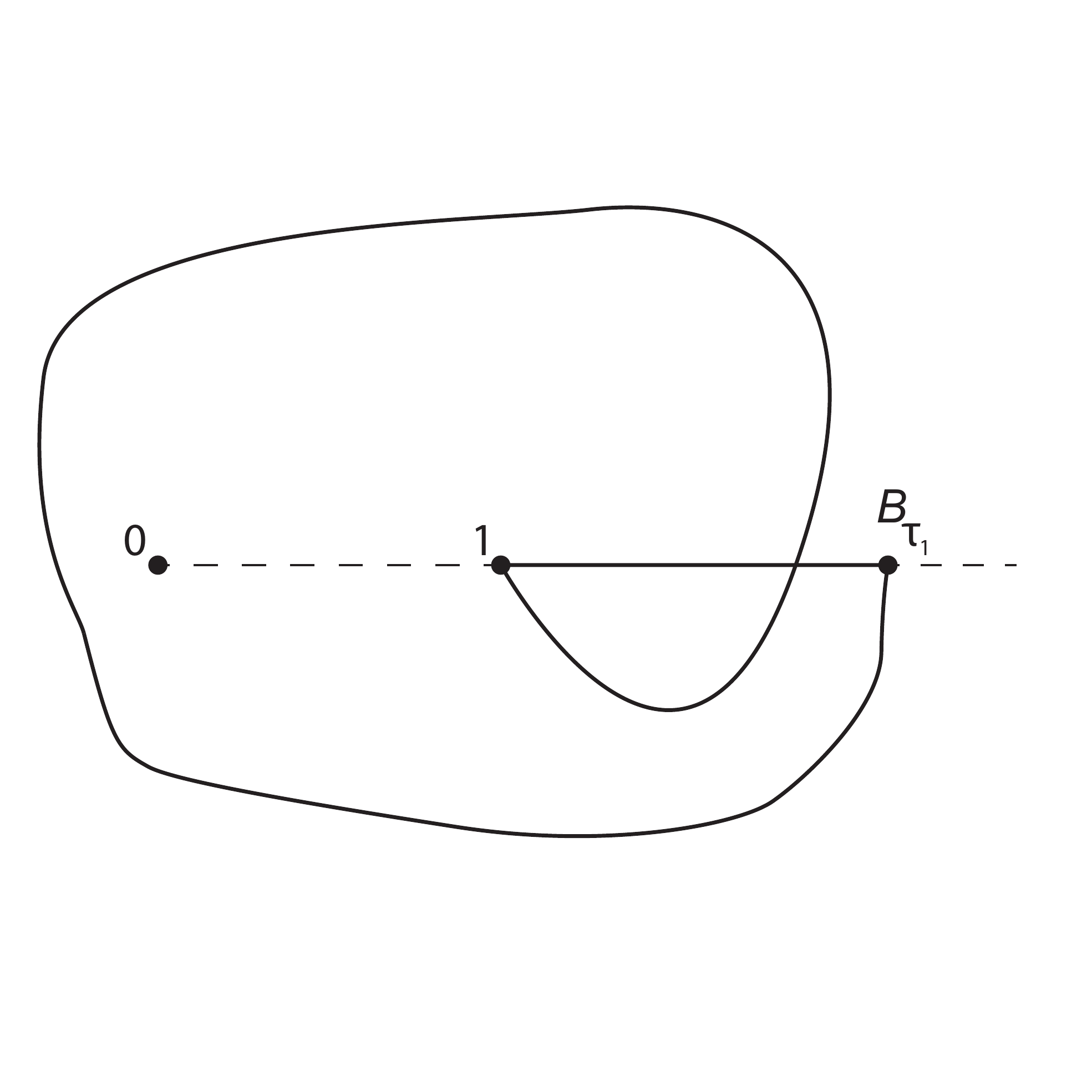}
\vspace{-.7in}
\caption{}
\end{wrapfigure}

\vski

The picture to the right shows the closed curve constructed in this manner from the example Brownian path shown in Figure 1. Since we have arrived at a contradiction, we conclude that our curve does intersect $\Om^c$ at some point, but since it does not do so on $\RR^+$ we see that it must have done so on the path $\{B_t:0 \leq t < \tau_1\}$, and this shows that $T_\Om < \tau_1$. Since $G_{\tau_1}(1,w) < \ff$ for any $w \neq 1$, monotonicity shows that $G_{\Om}(a,w) < \ff$ for all $w \neq a$ in $\Om$. Having obtained existence, we can obtain the harmonicity of $G_{\Om}(a,w)$ by noting that $\rho_t^\Om(a,w)$ evolves according to the heat equation $\dd_t \rho_t^\Om = \frac{1}{2} \Delta \rho_t^\Om$, so that

\begin{equation} \label{}
\begin{split}
\Delta G_{\Om}(a,w) & = \Delta \int_0^\ff \rho_t^\Om(a,w) dt = \int_0^\ff \Delta \rho_t^\Om(a,w) dt \\
& = 2 \int_0^\ff \dd_t \rho_t^\Om(a,w) dt = 2(\rho_\ff^\Om(a,w)- \rho_0^\Om(a,w)),
\end{split}
\end{equation}

where $\rho_\ff^\Om(a,w) = \lim_{t \lar \ff}\rho_t^\Om(a,w)$. But $\rho_\ff^\Om(a,w)$ and $\rho_0^\Om(a,w)$ by monotonicity are no more than $\rho_\ff^\CC(a,w)$ and $\rho_0^\CC(a,w)$, and both of these quantities are zero, so it follows that $\Delta G_{\Om}(a,w) = 0$, and $G_\Om$ is harmonic (The reader unsatisfied with the heuristic nature of this calculation is referred to \cite{mortper}, where harmonicity of $G_\Om$ is proved more carefully). To show that $G_{\Om}(a,w)$ has the correct logarithmic singularity at $w = a$, note that $G_{\Om}(a,w)$ is bounded above by $G_{\tau_1}$, modulo the M\"obius transformation, and below by $G_{D(a,\eps)}$ whenever $\epsilon>0$ is chosen so that $D(a, \eps) \subseteq \Om$. $G_{\tau_1}$ was shown in the previous section to have the correct order singularity at $w=a$ (the M\"obius transformation does not change this), and it may be checked that $G_{D(a,\eps)}$ does as well by applying the conformal map $z \lar a + \eps z$ to $G_{\DD}$ as calculated in the previous section. We see that $G_\Om(a,w) + \frac{1}{\pi} \ln |a-w|$ is bounded above and below in a neighborhood of $a$, and this implies that it may be extended to be harmonic at $a$ (see for example \cite[p. 166]{ahl}).

\vski

It remains only to show that $G_{\Om}(a,w) \lar 0$ as $w \lar \dd \Om$. In fact, the argument that $G_{\Om}(a,w)<\ff$ gives us a hint as to how to proceed, as we saw in the previous section that $G_{\tau_1}(1,w) \lar 0$ whenever $w \lar 0, \ff$. Thus, we have already shown that $G_{\Om}(a,w) \lar 0$ whenever $w \lar p$ if $p$ is a boundary point such that the line segment connecting $p_1$ to $a$ lies completely within $\Om$ (except at $p$). Let us now extend this argument to show that $G_{\Om}(a,w) \lar 0$ as $w \lar p$ for any $p \in \dd \Om$. We will begin by making the following simplifying assumption: let us assume that for any $p \in \dd \Om$ there is a curve $C$ connecting $a$ to $p$ which lies entirely within $\Om$ (except at $p$) such that $arg(w-p)$ remains bounded for all $w \in C$, where the $arg$ function is defined by continuity on the curve. In other words, we assume that any boundary point $p$ can be reached by a curve lying in $\Om$ and which does not wind infinitely many times around $p$. It is clear that this assumption would be satisfied by most simply connected domains that one would be likely to encounter in practice, although it fails for instance when $\Om$ is the complement of a logarithmic spiral. We will now argue that there is some positive integer $n$ such that $T_\Om < \tau_n$.

\vski

Extend $C$ by connecting $a$ to another boundary point $p_2 \neq p$ in any manner so that $arg(w-p)$ still remains bounded for all $w \in C$; this can always be done, by for instance drawing a straight path from $a$ in the direction opposite $p$ (assuming $p \neq \ff$), and letting $p_2$ be the first point encountered in $\dd \Om$ (possibly $\ff$). Again applying a M\"obius change of coordinates allows us to assume $a=1, p=0, p_2=\ff$, and $C$ is now a curve traveling from $0$ to $\ff$, passing through $1$, lying entirely within $\Om$ (except at $0, \ff \in \dd \Om$), such that $arg(w)$ remains bounded for all $w \in C$. These properties imply that $C$ can be realized as the image under the exponential map $z \to e^z$ of a curve $\hat C$ in $\CC$ such that $Re(w)$ ranges from $-\ff$ to $+\ff$ but $Im(w)$ remains bounded above and below as $w$ ranges over all points in $\hat C$. Furthermore we may translate $\hat C$ by any integer multiple of $2\pi i$ due the periodicity of the exponential function, and $\hat C$ must pass through $2\pi i m$ for some integer $m$. Thus, we may choose a positive integer $m_+$ and a realization $\hat C_+$ of $\hat C$ such that $\hat C_+$ passes through $2\pi i m_+$ but does not intersect $\RR$. Similarly, we may choose a negative integer $m_-$ and a realization $\hat C_-$ of $\hat C$ such that $\hat C_-$ passes through $2\pi i m_-$ but does not intersect $\RR$. Having chosen these, we may choose a positive integer $n$ such that the curves $\hat C_+, \hat C_-$ are entirely contained within $\{-2\pi n < Im(w) <2 \pi n\}$. It follows that the curves $\hat C_+, \hat C_-$ separate $0$ from the lines $\{Im(w) = \pm 2 \pi n\}$, and therefore if we start a Brownian motion $\hat B$ at $0$ and set $\hat \tau_n = \inf\{t \geq 0 : Im(B_t) = \pm 2 \pi n\}$, we see that we are guaranteed that $\hat B$ hits at least one of $\hat C_+, \hat C_-$ at a point $v$ before $\hat \tau_n$. This setup is illustrated in Figure 3.

\begin{center}
\begin{figure}
\vspace{-.7in}
\includegraphics[width=5 in,height=5 in]{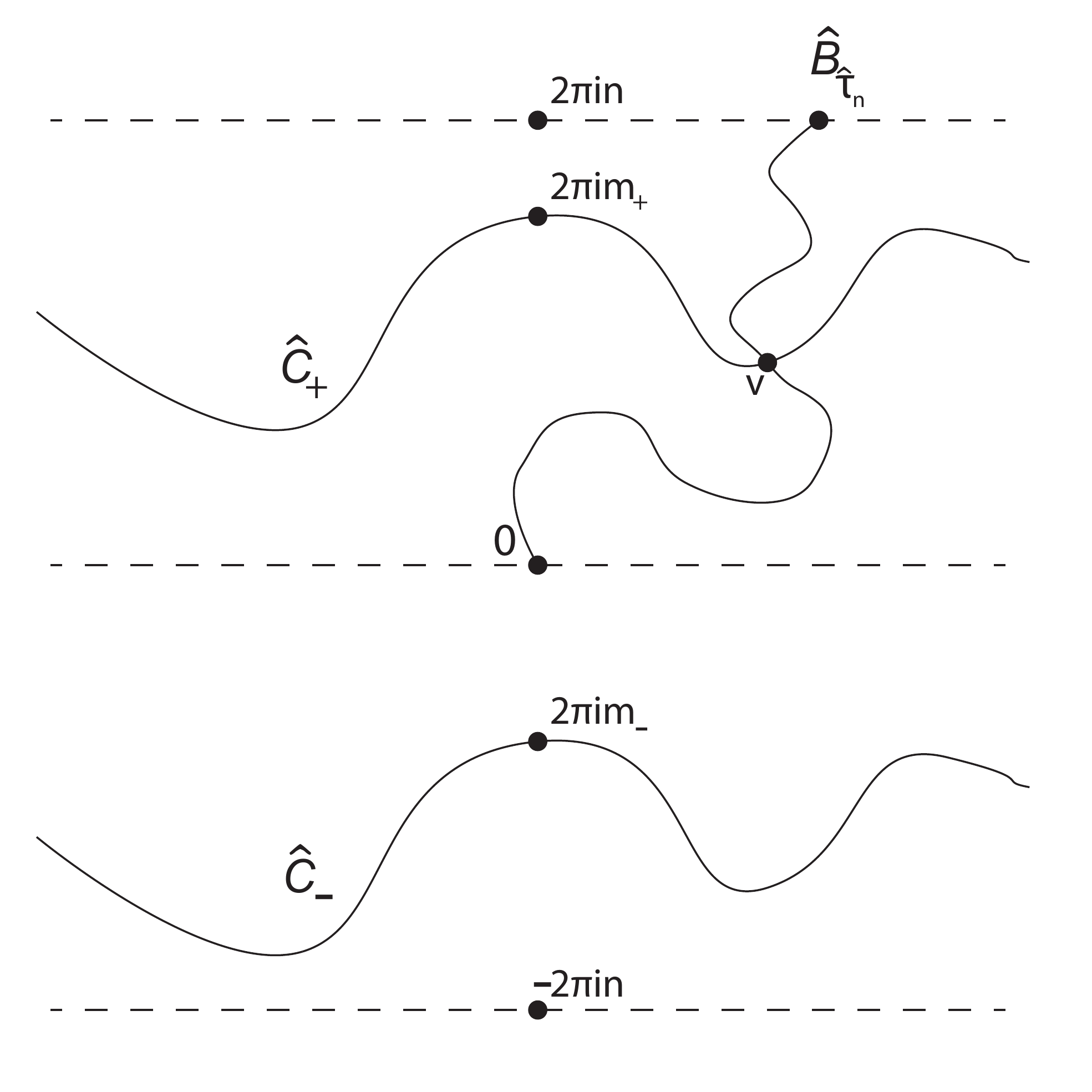}
\vspace{-.3in}
\caption{}
\end{figure}
\end{center}

Let us assume that $v \in \hat C_+$. We may now form a curve by following the Brownian path from $0$ to $v$, and then following $\hat C_+$ to $2\pi i m_+$; the image of this curve under the exponential map will be a closed curve in $\CC \bsh \{0\}$ which is not homotopic to a point. We conclude again that this curve must intersect $\Om^c$ at some point, but since it does not on $C$, which is the image of $\hat C_+$ under the exponential map, it must do so on the image of the Brownian path. We see that our Brownian motion starting at $a$ must exit $\Om$ before reaching the stopping time which is the projection of $\hat \tau_n$ under the exponential map; but this is simply the stopping time $\tau_n$. Thus, $T_\Om < \tau_n$. Since $G_{\tau_n}(1,w) \lar 0$ as $w \lar 0$, we obtain $G_{\Om}(a,w) \lar 0$ as $w \lar p$.

\vski

It now remains only to remove the assumption that any boundary point can be reached without winding infinitely many times around it. For positive integer $m$, let $\Om_m$ be the component of $\{w \in \Om: |w|<m, |p-w|>\frac{1}{m} \mbox{ for all }p \in \Om^c\}$ which contains $a$. It is easy to see that $\Om_m$ is simply connected and satisfies our simplifying assumption (because any point in $\dd \Om_m$ lies in $\Om$), so we have shown that there is a conformal map $f_m$ from $\DD \lar \Om_m$, taking $0$ to $a$. A standard normal families and Hurwitz's theorem argument shows that we have a function $f$ which is the locally uniform limit of a subsequence of the $f_m$'s, and that $f$ maps $\DD$ conformally onto $\Om$. This competes the proof.

\section{Infinite product identities}

In this section we show how calculating the Green's functions for certain domains in two different ways can yield identities, specifically several infinite products for trigonometric functions. Let us first map $\HH$ onto $\DD^\times = \DD \bsh \{0\}$ with the covering map $f(z) = e^{iz}$. Let $z=ai, w=b+ci$. Applying Proposition \ref{chimass} and using the calculation from Section 2 for $G_\HH(z,w)$ we get

\begin{equation} \label{}
\begin{split}
G_{\DD^\times}(e^{-a},e^{-c+bi}) & = \frac{1}{\pi} \sum_{n=-\ff}^{\ff} \ln \frac{|(a+c)i - (b+2\pi n)|}{|(a-c)i - (b+2\pi n)|} \\
& = \frac{1}{2\pi} \sum_{n=-\ff}^{\ff} \ln \frac{(b+2\pi n)^2 + (a+c)^2}{(b+2\pi n)^2 + (a-c)^2}.
\end{split}
\end{equation}

However, planar Brownian motion does not see points, i.e. $P_w(B_t = 0 \mbox{ for some }t \geq 0) = 0$ for any $w \neq 0$. Thus, this must agree with our earlier calculation that $G_{\DD}(e^{-a},e^{-c+bi}) = \frac{1}{2\pi} \ln \Big|\frac{1-e^{-a-c+bi}}{e^{-a} - e^{-c+bi}}\Big|^2$. Multiplying each by $2\pi$ and exponentiating gives the identity

\begin{equation} \label{mirror}
\prod_{n=-\ff}^\ff \frac{(b+2\pi n)^2+(a+c)^2}{(b+2\pi n)^2+(a-c)^2} = \Big|\frac{1-e^{-a-c+bi}}{e^{-a} - e^{-c+bi}}\Big|^2
\end{equation}

Special cases can take more familiar forms. If we set $b=0$ and rearrange a bit, we are led immediately to

\begin{equation} \label{}
\frac{(\frac{a+c}{2})^2}{(\frac{a-c}{2})^2} \Big(\prod_{n=1}^\ff \frac{(1+(\frac{a+c}{2\pi n})^2)}{(1+(\frac{a-c}{2\pi n})^2)}\Big)^2 = \frac{\sinh (\frac{a+c}{2})^2}{\sinh (\frac{a-c}{2})^2}.
\end{equation}

Multiplying both sides by $(\frac{a-c}{2})^2$ and taking the limit as $c \lar a$ yields the infinite product representation for $\sinh$:

\begin{equation} \label{}
\sinh a = a \prod_{n=1}^\ff (1+(\frac{a}{\pi n})^2).
\end{equation}

Returning to \rrr{mirror}, take now $b=\pi, c=a$, and we may reduce easily to

\begin{equation} \label{}
\cosh a = \prod_{n=1}^\ff (1+(\frac{a}{\pi (n-1/2)})^2),
\end{equation}

Note that the infinite product representations for sine and cosine can be derived from these, as $\sin a = i \sinh(-ia)$ and $\cos a = \cosh (ia)$, and we obtain:

\begin{equation} \label{}
\sin a = a \prod_{n=1}^\ff (1-(\frac{a}{\pi n})^2); \qquad \cos a = \prod_{n=1}^\ff (1-(\frac{a}{\pi (n-1/2)})^2).
\end{equation}

\vski

Now let us consider the strip $U= \{-1<Im(z)< 1\}$. Arguing similarly to the upper half-plane case, $\rho_t^U (z,w)$ must capture the probability of Brownian paths near $w$ at time $t$, but only those which have not yet touched $\{Im(z) = -1\}$. Applying the reflection principle as before we see that we must subtract $\rho_t^\CC (z,-2i + \bar w)$, since $-2i + \bar w$ is the reflection of $w$ over $\{Im(z) = -1\}$. We must also get rid of paths which touch $\{Im(z) = 1\}$ before proceeding to $w$, so we need to subtract $\rho_t^\CC (z,2i + \bar w)$, since $2i + \bar w$ is the reflection of $w$ over $\{Im(z) = 1\}$; however in doing so we have subtracted too much, since we have two times subtracted the probability corresponding to paths which touch both of $\{Im(z) = -1\}$ and $\{Im(z) = 1\}$. In order to rectify this, we must add in the probability corresponding to such paths, and applying the reflection principle again we see that this can be done by reflecting $w$ twice in both directions (once over $\{Im(z) = \pm 1\}$, and then once over $\{Im(z) = \pm 3\}$), i.e. adding $\rho_t^\CC (z,-4i + w)$ and $\rho_t^\CC (z,4i + w)$. But now we've added too much, again twice counting paths that travel between $\{Im(z) = \pm 1\}$ three times, so we must reflect again and subtract, etc. Continuing in this manner, we conclude that

\begin{equation} \label{wavrdmas}
\begin{split}
\rho_t^U (z,w) & = \rho_t^\CC (z,w) - (\rho_t^\CC (z,-2i + \bar w) + \rho_t^\CC (z,2i + \bar w)) \\
& \quad + (\rho_t^\CC (z,4i + w) + \rho_t^\CC (z,-4i + w)) - (\rho_t^\CC (z,6i + \bar w) + \rho_t^\CC (z,-6i + \bar w)) + \ldots
\end{split}
\end{equation}
Integrating, using the same method as in the half-plane case, we obtain

\begin{equation} \label{}
\begin{split}
G_{U}(z,w) & = \frac{1}{\pi} \sum_{n=-\ff}^\ff \ln \Big( \frac{|z -((4n+2)i + \bar w)|}{|z-(4ni +w)|}\Big)\\
& = \frac{1}{\pi} \sum_{n=-\ff}^\ff \ln \Big( \frac{|z -((4n-2)i + \bar w)|}{|z-(4ni +w)|}\Big).
\end{split}
\end{equation}

If $Re(z-w) = a, Im(z)=b, Im(w)=c$, then this is

\begin{equation} \label{sonya}
G_{U}(z,w)  = \frac{1}{2\pi} \sum_{n=-\ff}^\ff \ln \Big( \frac{a^2 + ((b+c)-(4n-2))^2}{a^2 + ((b-c)-4n)^2}\Big).
\end{equation}

We can also calculate $G_{U}(z,w)$ by mapping $U$ conformally to the unit disk by $f(z)=\tan(\frac{\pi i}{4}z)$ (alternatively, one may use the map $f(z) = ie^{\frac{\pi}{2}z}$, which maps $U$ conformally to $\HH$, in conjunction with \rrr{hp}). Using \rrr{lamb} we obtain

\begin{equation} \label{}
G_{U}(z,w)  = \frac{1}{\pi} \ln \Big| \frac{1 + \tan(\frac{\pi b}{4})\tan(\frac{\pi i}{4}a - \frac{\pi}{4}c)}{\tan(\frac{\pi b}{4}) + \tan(\frac{\pi i}{4}a - \frac{\pi}{4}c)}\Big|.
\end{equation}

Equating these two expressions and exponentiating leads to

\begin{equation} \label{4hands}
\prod_{n=-\ff}^\ff \Big( \frac{a^2 + ((b+c)-(4n-2))^2}{a^2 + ((b-c)-4n)^2}\Big) = \Big| \frac{1 + \tan(\frac{\pi b}{4})\tan(\frac{\pi i}{4}a - \frac{\pi}{4}c)}{\tan(\frac{\pi b}{4}) + \tan(\frac{\pi i}{4}a - \frac{\pi}{4}c)}\Big|^2.
\end{equation}

Special cases take on simpler forms. For instance, setting $b=c=0$ and replacing $a$ with $\frac{-4ia}{\pi}$ leads to the identity

\begin{equation} \label{lynne}
\tan^2a = -\prod_{n=-\ff}^\ff \frac{n^2 - \frac{a^2}{\pi^2}}{(n-1/2)^2 - \frac{a^2}{\pi^2}}.
\end{equation}

On the other hand, setting $a=b=0$ and replacing $c$ with $\frac{-4c}{\pi}$ leads to the identity

\begin{equation} \label{sonya2}
\tan c = -\prod_{n=-\ff}^\ff \frac{n + \frac{c}{\pi}}{(n-1/2) - \frac{c}{\pi}}.
\end{equation}

{\bf Remarks:} The reader unhappy with the heuristic derivation of \rrr{wavrdmas} may prefer to argue as follows. Let us define $\tau(b) = \inf\{t \geq 0: Im(B_t) = b\}$, and then recursively define $\tau(b_1, b_2, \ldots, b_{n+1}) = \inf\{t \geq \tau(b_1, \ldots, b_{n}): Im(B_t) = b_{n+1}\}$; that is, $\tau(b_1, b_2, \ldots, b_{n})$ is the first time at which $Im(B_t)$ has visited the sequence $b_1, b_2, \ldots , b_n$ in order. We then have almost surely, for any $A \subseteq U$,

\begin{equation} \label{jess}
\begin{split} \nonumber
1_{\{B_t \in A, t < T_U\}} & = 1_{\{B_t \in A, t < \tau(1) \wedge \tau(-1)\}} \\
& = 1_{\{B_t \in A\}} - (1_{\{B_t \in A, t > \tau(1)\}} + 1_{\{B_t \in A, t > \tau(-1)\}}) \\
& \quad + (1_{\{B_t \in A, t > \tau(1,-1)\}} + 1_{\{B_t \in A, t > \tau(-1,1)\}}) - (1_{\{B_t \in A, t > \tau(1,-1,1)\}} + 1_{\{B_t \in A, t > \tau(-1,1,-1)\}}) + \ldots
\end{split}
\end{equation}

Note that the right hand side of \rrr{jess} has only a finite number of nonzero terms for fixed $t$ and Brownian path. \rrr{wavrdmas} then follows by letting $A$ be a small disk centered at $w$, diving by the area of $A$, taking expectations, and letting the radius of $A$ shrink to 0. 

\vski

It should also be mentioned that an analytic method equivalent to the method of reflection used in the infinite strip example above, including the derivation of an infinite product equivalent to \rrr{4hands}, has previously appeared in \cite{melgreen} and \cite{melgreen2}.

\section{Concluding remarks}

The proof of the Riemann mapping theorem given here is not the first to make use of the winding of Brownian motion: \cite{bass} contains one as well, although the details are fairly different. The proof given there uses the support theorem (which states that Brownian motion uniformly approximates any path on a finite time interval with positive property), which in turn depends on a Girsanov transformation of Wiener measure, in order to prove the regularity of the boundary of any simply connected domain, thence to the existence of the Green's function. In contrast, the proof given here does not make use of stochastic calculus (other than implicitly in the application of L\'evy's theorem).

\vski

The Riemann mapping theorem has an interesting history, as is detailed for instance in \cite{walshhistory}. Riemann formulated the theorem (with smoothness assumptions on the boundary) and gave a proof in 1851, although his proof is now considered to have been incomplete at the time. Following this, such noteworthy mathematicians as Weierstrass, Hilbert, Schwarz, Harnack, and Poincar\'e worked to correct the proof and relax the boundary restrictions, until it was finally proved in full generality by Osgood in 1900 \cite{osgood}. As the name of his paper suggests, Osgood's proof (as well as the earlier partial ones) proceeded by showing the existence of the analyst's Green's function for any simply connected domain, and thence to the conformal map. This is opposed to the standard modern proof, developed and simplified over the years by such mathematicians as Carath\'eodory, Koebe, Fej\'er, Riesz, and Ostrowski, which depends on function theoretic arguments and bypasses Green's function. Clearly, the proof given here is more along the lines of the early proofs, and it is interesting to note that regularity conditions of the boundary again come into play, although the ease with which boundary points may be approached from within the domain was the issue for us, rather than any sort of smoothness.

\vski

Although the probabilist's Green's function is known to be harmonic whenever the stopping time in question is the exit time of a domain (see for instance \cite{mortper}), this does not hold for arbitrary stopping times. A good example of this is for the stopping time $\tau_1$: it can be shown that the support of $G_{\tau_1}(1,w)$ is $\CC \bsh \{0\}$ but $G_{\tau_1}(1,w)$ is not harmonic for $w \in \RR^+$. The reason for this is that the Brownian motion is stopped on $\RR^+$. In general, if the stopping time $\tau$ is the projection under an analytic function of an exit time, then the Green's function of that stopping time will be harmonic at all points in its support other than those in the support of $\{B_\tau\}$.

\vski

A possible application of the method of proof given here may be in acquiring bounds and rates of decay for the analyst's Green's function. For instance, if the domain in question is starlike with $a$ in the center (any point in the domain can be connected to $a$ with a straight line segment contained in the domain), then the Green's function for the domain can be bounded by $G_{\tau_1}(1,w)$ (which can be calculated explicitly) composed with a suitable M\"obius transformation. Details are omitted, but should be straightforward to supply if such a result is required.

\vski 

It may be tempting to search for new domains, stopping times, and analytic functions to which to apply Proposition \ref{chimass} in order to find new infinite product representations. The reader determined to undertake this should be warned that conformal invariance means that many seemingly unrelated constructions are actually equivalent. A good example of this is that the function $\sin z$ maps the upper half-plane $\HH$ onto $\CC \backslash [-1,1]$, essentially wrapping the half-plane in a periodic manner about the line segment. Since the Green's function for $\CC \backslash [-1,1]$ can be calculated directly, it seems as though we will obtain a new identity, however we obtain simply \rrr{mirror}, which was obtained in the punctured disk example above. This is because $\CC \backslash [-1,1]$ can be conformally mapped to the punctured disk, with the point at $\ff$ corresponding to $0$, and the function $\sin z$ essentially keeps track of the homotopy classes of curves at $w$ at time $t$, exactly as the exponential function does in the punctured disk. For another example, if we start our Brownian motion at $1$ with $arg(B_t)=0$, and let $\tau_K = \inf\{t \geq 0 : |arg(B_t)| = K\}$, then using reflection we may express $\rho^{\tau_K}_t(1,w)$ as an infinite sum of terms of the form $\rho^{\hat \tau_{\hat K}}_t(1,w)$, where $\hat \tau_{\hat K} = \inf\{t \geq 0 : arg(B_t) = \hat K\}$. However, this is simply a transformation by the exponential function of the argument used in the infinite strip above, and therefore yields the same infinite product identity.

\section{Acknowledgements} I would like to thank Andrea Collevecchio, Burgess Davis, and Michael Kozdron for helpful conversations, and Jiro Akahori for a helpful invitation. It should be mentioned that the proof of the Riemann mapping theorem given here is inspired by, and similar in spirit to, the proof of Picard's theorem by Davis (see \cite{davispicard}), which also makes use of the winding of Brownian motion. I am also grateful for support from Australian Research Council Grants DP0988483 and DE140101201.

\bibliography{CAbib}

\end{document}